  \documentclass[12pt]{article}

\usepackage{amscd}
\usepackage{amsmath}
\usepackage{amssymb}
\usepackage{amsthm}
\usepackage{epsf}
\usepackage{geometry}
\usepackage{latexsym}
\usepackage{verbatim}
\usepackage{graphicx}
\usepackage{amssymb}
\usepackage{epstopdf}
\input epsf.tex

\newtheorem{theorem}{Theorem}
\newtheorem{lemma}[theorem]{Lemma}
\newtheorem{definition}[theorem]{Definition}

\newtheorem{remark}[theorem]{Remark}

\newcommand{\re}{\operatorname{Re} \ }
\newcommand{\im}{\operatorname{Im} \ }
\newcommand{\n}{\noindent}

\begin{document} 

\title{The Weierstrass Representation always gives a minimal surface}

\bibliographystyle{plain}

\author{Roshan Sharma
 \thanks {Williams College, Williamstown, MA 01267 USA.
    roshan.sharma@williams..edu}}

\maketitle 

\begin{abstract}
We give a simple, direct proof of the easy fact about the Weierstrass Representation, namely, that it always gives a minimal surface. Most presentations include the much harder converse that every simply connected minimal surface is given by the Weierstrass Representation.
\end{abstract}
\section{Introduction}
The celebrated Weierstrass Representation (Theorem 3) gives all simply connected minimal surfaces in terms of two arbitrary complex functions (sometimes allowing poles). Here we give a simple direct proof of the easy direction: that every surface given by the formula is a minimal surface. Minimal surfaces by definition satisfy the variational condition for minimizing area, namely that the mean curvature vanishes. Nice examples are soap films (Figures 1 and 2), delicately balanced with their principal curvatures equal and opposite.  
\clearpage
\begin{figure}
\begin{center}\includegraphics[scale = 0.75]{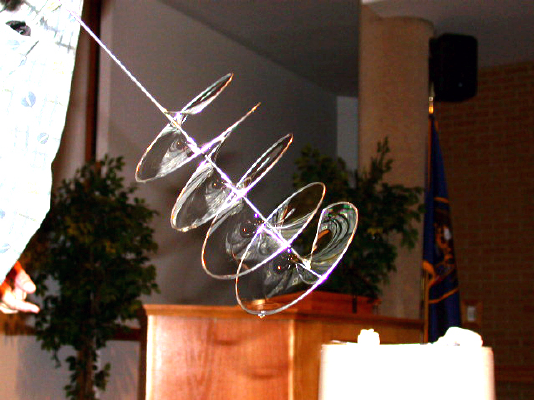}\end{center}
\scriptsize{\caption{A helicoid formed by a soap film for the wire boundary. (academic.csuohio.edu/oprea\textunderscore j/utah/Prospects.html Accessed 8/15/12. Used by permission, all rights reserved.)}}
\end{figure}
\medskip
\begin{figure}
\begin{center}\includegraphics[scale = 0.75]{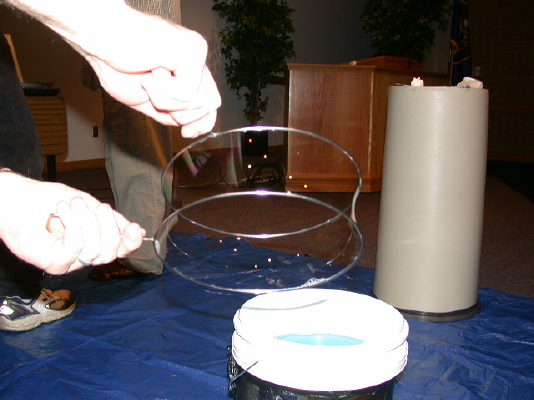}\end{center}
\scriptsize{\caption{A catenoid formed by a soap film for the wire boundary. (academic.csuohio.edu/oprea\textunderscore j/utah/Prospects.html Accessed 8/15/12. Used by permission, all rights reserved.)}}
\end{figure}
\subsection{Acknowledgments}
This paper is based on my Williams senior colloquium, advised by Professor Frank Morgan.
\section{The Weierstrass Representation}
We start with some definitions.
\begin{definition} \textup{A {\it {minimal surface}} in $\mathbb{R}^3$ has mean curvature 0 at every point. \\ \\ The {\it{square}} of a complex vector ${\bf v} = (v_1, v_2, v_3)$ is defined as }
\begin{eqnarray*}
{\bf v}^2 = {\bf v}\cdot {\bf v} = v_1^2 + v_2^2 + v_3^2 .
\end{eqnarray*}
\end{definition}
\n The following lemma gives the key property of the Weierstrass Representation.
\begin{lemma}
For complex numbers or functions $\phi  = \left[\begin{matrix}f(1-g^2) \\ if(1+g^2) \\ 2fg\end{matrix}\right]$, 
the following statements are true. \begin{description}
\item[(1)]$(\re \phi)^2 - (\im \phi)^2=0$.
\item[(2)] $(\re \phi) \cdot (\im \phi) = 0$.
\end{description}
\begin{proof}
First we note that
\begin{eqnarray*}
\phi^2 &=& f^2 - 2f^2g^2 + f^2g^4 - f^2 - 2f^2g^2 - f^2g^4 + 4f^2g^2 \\
&=& 0.
\end{eqnarray*}
Therefore, we get \begin{eqnarray*}
0 = \re \phi^2 = (\re \phi)^2 - (\im \phi)^2,
\end{eqnarray*}
because $\re (a+ib)^2 = a^2 - b^2$. Similarly, 
\begin{eqnarray*}
0 = \im \phi^2  = 2(\re \phi)\cdot (\im \phi),
\end{eqnarray*}
because $\im (a+ib)^2 = 2ab$.
\end{proof}
\end{lemma}
\n We now give our main result.
\begin{theorem}
For any complex functions $f(z)$ and $g(z)$ on the unit disk or complex plane, the surface ${\bf x}(z)$ is minimal, where ${\bf x}$ is the real part of an integral of 
\begin{eqnarray*}
\phi =  \left[\begin{matrix} f(1-g^2) \\ if(1+g^2) \\ 2fg\end{matrix}\right].
\end{eqnarray*} 
\end{theorem}
\begin{remark} 
Actually you can allow $g(z)$ to have poles as long as $fg^2$ is analytic. 
\end{remark}
\begin{proof}[Proof of theorem]
For a surface ${\bf x}(z)$ in $\mathbb{R}^3$, the mean curvature $H$ is given by \cite{book}
\begin{eqnarray}
H = P_\perp \dfrac{{\bf x}_v^2 {\bf x}_{uu} - 2({\bf x}_u\cdot {\bf x}_v){\bf x}_{vv} + {\bf x}_u^2 {\bf x}_{vv}}{{\bf x}_u^2{\bf x}_v^2 - ({\bf x}_u\cdot {\bf x}_v)^2},
\end{eqnarray}
where $P_\perp$ denotes projection onto the line normal to the surface and the subscripts $u, v$ denote partial differentiation with respect to the real and imaginary parts of $z$. Letting $\Phi$ denote an integral of $\phi$, we compute that 
\begin{eqnarray}
{\bf x}_u = \re [\Phi_u] = \re\left[\dfrac{d \Phi}{dz} \dfrac{\partial z}{\partial u}\right] = \re \phi,
\end{eqnarray} as $z = u+iv$ so $\dfrac{\partial z}{\partial u} = 1$. Similarly, 
\begin{eqnarray}
{\bf x}_v = \re [\Phi_v] = \re\left[\dfrac{d \Phi}{dz}\dfrac{\partial z}{\partial v}\right] = \re(i\phi) = - \im \phi, 
\end{eqnarray}
because $\re [i(a+ib)] = -b$. By Lemma 1(2), we get
\begin{eqnarray}
{\bf x}_u \cdot {\bf x}_v = 0.
\end{eqnarray}
Since
\begin{eqnarray}
{\bf x}_{u}^2 &=&  (\re \phi)^2
\end{eqnarray}
and
\begin{eqnarray}
{\bf x}_{v}^2 &=& (\im \phi)^2,
\end{eqnarray}
by Lemma 1(1), we get 
\begin{eqnarray}
{\bf x}_{u}^2 = {\bf x}_{v}^2.
\end{eqnarray}
Furthermore, since
\begin{eqnarray}
{\bf x}_{uu} = \dfrac{d {\bf x}_u}{du} = \dfrac{d {\bf x}_u}{dz}\dfrac{\partial z}{\partial u}  = \re (\phi')
\end{eqnarray}
and 
\begin{eqnarray}
{\bf x}_{vv} = \dfrac{d {\bf x}_v}{dv} = \dfrac{d {\bf x}_v}{dz}\dfrac{\partial z}{\partial v} = -\re (\phi'),
\end{eqnarray}
therefore,
\begin{eqnarray}
{\bf x}_{uu} + {\bf x}_{vv} = 0.
\end{eqnarray}
Plugging Equations (4, 7, 10) into formula (1) yields $H = 0$, the definition of a minimal surface.
\end{proof}
\begin{remark}
\textup{The converse of Theorem 3 also holds: every simply connected minimal surface is given by the Weierstrass Representation [1, Theorem 18]}. 
\end{remark}

\end{document}